\input amstex
\loadeufm

\documentstyle{amsppt}

\magnification=\magstep1

\baselineskip=20pt
\parskip=5.5pt
\hsize=6.5truein
\vsize=9truein
\NoBlackBoxes

\define\br{{\Bbb R}}
\define\oo{{\omega}}
\define\e{{\varepsilon}}
\define\OO{{\Omega}}
\define\CL{{\Cal{L}}}
\topmatter
\title
Bounds of Riesz Transforms on  $L^p$ Spaces
for Second Order Elliptic Operators
\endtitle

\author Zhongwei Shen
\endauthor

\leftheadtext{ZHONGWEI SHEN}
\rightheadtext{Bounds of Riesz Transforms on $L^p$ Spaces}

\address Department of Mathematics, University of Kentucky,
Lexington, KY 40506.
\endaddress

\email shenz\@ms.uky.edu
\endemail


\abstract
Let $\CL=-\text{div}(A(x)\nabla)$ be a second order elliptic
operator with real, symmetric, bounded measurable
 coefficients on $\br^n$ or on a bounded Lipschitz domain
subject to Dirichlet boundary condition. For any fixed $p>2$,
a necessary and sufficient
condition is obtained for the boundedness
of the Riesz transform $\nabla (\CL)^{-1/2}$ on the $L^p$ space.
 As an application, for $1<p<3+\e$, we  establish 
the $L^p$ boundedness of Riesz transforms on Lipschitz domains
for operators with $VMO$ coefficients.
The range of $p$ is sharp. 
The closely related boundedness of $\nabla (\CL)^{-1/2}$
on weighted $L^2$ spaces is also studied.
\endabstract

\endtopmatter

\document

\centerline{\bf 1. Introduction}

Consider the second order elliptic operator of divergence form
$$
\CL=-\text{div}\big(A(x)\nabla\big)\ \ \ \ \text{ on } \ \ \OO,
\tag 1.1
$$
where $\OO=\br^n$ or a bounded open set of $\br^n$.
In the case of bounded domains, we impose the Dirichlet boundary
condition $u=0$ on $\partial\OO$.
Throughout this paper, we assume that
 $A(x)=(a_{jk}(x))$ is an $n\times n$ symmetric matrix
with real-valued, bounded measurable entries
satisfying the uniform ellipticity condition,
$$
\mu |\xi|^2
\le \sum_{j,k=1}^n a_{jk} (x) \xi_j \xi_k 
\le \frac{1}{\mu} |\xi|^2,
 \text{ for all }\xi,\, x \in \br^n
\ \text{and some }\mu>0.
\tag 1.2
$$
Under these assumptions, it is known that the Riesz transform
$\nabla (\CL)^{-1/2}$ is bounded on $L^p(\OO)$ for $1<p<2+\e$
where $\e=\e(n,\mu)>0$, and is 
of weak type $(1,1)$ (see e.g. \cite{4,7}).
Moreover, the range of $p$ is sharp.
The main purpose of this paper is to investigate
 the $L^p$ boundedness of
the Riesz transform for $p>2$, as well as the closely related
boundedness  on weighted $L^2$ spaces, under some {\it additional}
 conditions.
For any fixed $p>2$, we obtain a necessary and sufficient condition
for the boundedness of the Riesz transform on $L^p(\OO)$.
Armed with this condition, for elliptic operators with
$VMO$  coefficients, we are able to establish the $L^p$ boundedness
of Riesz transforms on Lipschitz domains 
for the optimal range of $p$.

Theorems A, B and C below are the main results of the paper.

\proclaim{\bf Theorem A}
Let $\CL$ be a second-order uniform elliptic operator of 
divergence form with real, symmetric, bounded measurable coefficients
on $\br^n$, $n\ge 2$.
For any fixed $p>2$, the following statements are equivalent.
\item{i)} There exists a constant $C>0$ such that
for any ball $B=B(x_0,r)$ and any $W^{1,2}_{loc}$
weak solution of
$\CL u=0$ in $3B=B(x_0,3r)$, one has $|\nabla u|\in L^{p}(B)$ and
$$
\left(\frac{1}{|B|}\int_B |\nabla u |^{p}dx\right)^{1/p}
\le C\, \left(\frac{1}{|2B|}
\int_{2B} |\nabla u|^2 dx\right)^{1/2}.
\tag 1.3
$$
\item{ii)} There exists $\e>0$ such that
the  Riesz transform $\nabla (\CL)^{-1/2}$ is 
bounded on $L^q(\br^n, dx)$ for any $1<q<p+\e$.
\item{iii)} There exists $\delta>0$ such that if $\oo$
is an $A_s$ weight with $s=2(1-\frac{1}{p}) +\delta$, then
$$
\| \nabla (\CL)^{-1/2} f\|_{L^2(\br^n, \frac{dx}{\oo})}
\le C\, \| f\|_{L^2(\br^n, \frac{dx}{\oo})},
\tag 1.4
$$
where $C$ depends on the $A_s$ bound of $\oo$.

\noindent In particular, $\nabla (\CL)^{-1/2}$ is bounded on $L^p(\br^n, dx)$
if and only if condition (i) holds for the same $p$. Consequently,
the set of exponents $p\in (1,\infty)$
for which $\nabla(\CL)^{1/2}$ is bounded on $L^p(\br^n, dx)$
is an open interval $(1,p_0)$ with $2<p_0\le \infty$.
\endproclaim

We remark that condition (ii) follows directly from (iii) by
an extrapolation theorem, due to Rubio de Francia.
It is also not hard to see that condition (i) follows from (ii)
by a standard localization argument, since the $L^p$ boundedness of
the Riesz transform yields the $W^{1,p}$ estimates
for $\CL$. The rest of the proof of Theorem A,
 however, is much more involved.
To prove that condition (i) implies (ii), we use a new
and refined version
 of the celebrated Calder\'on-Zygmund
Lemma.
See Theorem 3.1. 
This theorem, 
formulated by the author in \cite{17}, was
inspired by a paper of Caffarelli and Peral \cite{6}
as well as a recent work of L.~Wang \cite{19}.
For any fixed $p>2$, it
gives a sufficient condition for an $L^2$ bounded sublinear operator
to be bounded on $L^q$ for all $2<q<p$.
It enables us to
show that the operator $\nabla \CL^{-1}\text{div}$
is bounded on $L^p$ under condition (i). The boundedness of
the Riesz transform then follows from the
fact that $\|\CL^{1/2} f\|_q\le C\, \|\nabla f\|_q$
for any $1<q<\infty$ \cite{4}.
To show that condition (i) implies (iii),
the basic observation is that condition (i) leads to an
$L^p$ estimate on the kernel function of the Riesz transform.
Using this estimate on the kernel as well as  the $L^p$
boundedness established above, we show that the sharp function of
the adjoint of the Riesz transform can be dominated
by the Hardy-Littlewood maximal function.
The desired estimate (1.4) then follows from the weighted
norm inequalities for the sharp and maximal functions.

It is not very 
difficult to extend the argument above to the case of bounded
Lipschitz domains. This gives us the following.

\proclaim{\bf Theorem B}
Let $\OO$ be a bounded Lipschitz domain in $\br^n$, $n\ge 2$. Let
$\CL$ be a second-order uniform elliptic operator of divergence
form on $\OO$, subject to Dirichlet boundary condition.
For any fixed $p>2$, the following statements are equivalent.
\item{i)}
There exist constants $C>1$, $\alpha_2>\alpha_1>1$ and $r_0>0$ 
such that for any ball $B(x_0,r)$
with the property that $0<r<r_0$ and either $x_0\in\partial\OO$
or $B(x_0,\alpha_2 r)\subset\OO$, and for any weak solution
of $\CL u =0$ in $\OO\cap B (x_0,\alpha_2 r)$ and $u=0 $ on 
$B(x_0, \alpha_2 r)\cap \partial\OO$ (if $x_0\in \partial\OO$),
one has $|\nabla u|\in L^{p}(\OO\cap B(x_0,r))$ and
$$
\left(\frac{1}{r^n}\int_{\OO\cap B(x_0,r)} |\nabla u|^{p} dx \right)^{1/p}
\le C\, \left(\frac{1}{r^n}
\int_{\OO\cap B(x_0,\alpha_1 r)} |\nabla u|^2 dx\right)^{1/2}.
\tag 1.5
$$
\item{ii)} There exists $\e>0$ such that
the Riesz transform $\nabla (\CL)^{-1/2}$ is
bounded on $L^q(\OO, dx)$ for any $1<q<p+\e$.
\item{iii)}
There exists $\delta>0$ such that if $\omega$
is an $A_s(\br^n)$ weight with $s=2(1-\frac{1}{p}) +\delta$, then
$$
\| \nabla (\CL)^{-1/2} f\|_{L^2(\OO, \frac{dx}{\oo})}
\le C\, \| f\|_{L^2(\OO, \frac{dx}{\oo})},
\tag 1.6
$$
where $C$ depends on the $A_s$ bound of $\oo$.

\noindent In particular, 
$\nabla (\CL)^{-1/2}$ is bounded on $L^p(\Omega, dx)$
if and only if condition (i) holds for the same $p$. Consequently,
the set of exponents $p\in (1,\infty)$
for which $\nabla(\CL)^{1/2}$ is bounded on $L^p(\Omega, dx)$
is an open interval $(1,p_0)$ with $2<p_0\le \infty$.
\endproclaim

A few remarks are in order.

\remark{\bf Remark 1.7} Let $\OO=\br^n$ or a bounded Lipschitz
domain.
It follows from Theorems A and B that $\nabla (\CL)^{-1/2}$
is bounded on $L^p(\OO, dx)$ for any $p\in (1,\infty)$ if and only
if it is bounded on $L^2(\OO, \omega dx)$ for any
$\oo \in A_2(\br^n)$.
To see this, we note that $\oo\in A_2$ implies that $\oo\in A_q$
for some $q<2$, and $\frac{1}{\omega}\in A_2$.
For the classical Riesz transform
$\nabla(-\Delta)^{-1/2}$,  the boundedness on  $L^p(\br^n)$
 for $1<p<\infty$, and on $L^2(\br^n, \oo dx)$
with $\oo\in A_2(\br^n)$ is well known (see e.g. \cite{18,10}).
\endremark

\remark{\bf Remark 1.8}
For a general second order elliptic operator $\CL$
with real, symmetric, bounded measurable coefficients,
$\nabla (\CL)^{-1/2}$ is bounded
on $L^p$ for $1<p<2+\e$.
The range of $p$ was shown to be optimal by
C.~Kenig (see \cite{4, pp.119-121}).
It is worth  mentioning that the
 boundedness of $\nabla (\CL)^{-1/2}$ on $L^p$
is equivalent to the inequality
$\| \nabla f\|_p\le C\, \| \CL^{1/2}f \|_p$.
The reverse inequality $\|\CL^{1/2} f\|_p
\le C\, \| \nabla f\|_p$, nevertheless, holds
for all $1<p<\infty$ \cite{4,5}.
The proof of Theorems A and B depends on this fact.
\endremark

\remark{\bf Remark 1.9} 
By a simple geometric
observation, one may see that condition (i) in Theorem B
is equivalent to the following. There exist $C_1>0$,
$\alpha_4>\alpha_3>1$ and $r_1>0$ such that
for any $D=B(x_0,r)\cap\OO\neq\emptyset$ with
$x_0\in \br^n$, $0<r<r_1$, and for any weak solution
of $\CL u=0$ in $ \OO\cap B(x_0,\alpha_4 r)$ and $u=0$
on $B(x_0,\alpha_4 r)\cap\partial\OO$ (if it's not empty), one has
$|\nabla u|\in L^p(D)$ and
$$
\left\{ \frac{1}{r^n}
\int_{D}
|\nabla u|^p \, dx\right\}^{1/p}
\le C_1\,
\left\{ \frac{1}{r^n}
\int_{\OO \cap B(x_0,\alpha_3 r)}
|\nabla u|^2\, dx\right\}^{1/2}.
\tag 1.10
$$
By the reverse H\"older inequality estimates
\cite{12, pp.122-123}, this implies that
condition (i) in Theorem B has the self-improvement property. This is, if
$\CL$ satisfies condition (i) in Theorem B for some
$p>2$, then it satisfies condition (i) for some $\bar{p}>p$.
Clearly, the same can be said about the condition (i)
in Theorem A.
It follows that the set of exponents $p\in (1,\infty)$
for which $\nabla (\CL)^{-1/2}$ is bounded on $L^p(\Omega)$
is an open interval.
\endremark

Let $\CL=-\Delta$ on a bounded Lipschitz domain 
$\OO$, subject to Dirichlet boundary condition.
Using the solvability of the $L^2$ regularity
problem and the boundary H\"older
estimates (see \cite{14}),
 it is not hard to show that condition (i) in Theorem
B holds for $p=3$ if $n\ge 3$, and for $p=4$ in the case $n=2$
(see Lemma 4.1). 
It follows that for $n\ge 3$,
the Riesz transform $\nabla (\CL)^{-1/2}$ is bounded on
$L^p(\OO)$ for $1<p<3+\e$, and on 
$L^2 (\OO,\frac{dx}{\omega})$ with
$\oo\in A_{\frac{4}{3}+\delta}(\br^n)$.
If $n=2$, $\nabla (\CL)^{-1/2}$ is bounded on $L^p(\OO)$
for $1<p<4+\e$, and on $L^2(\OO,\frac{dx}{\oo})$ with
$\oo\in A_{\frac{3}{2}+\delta}(\br^2)$.
The ranges of $p$ are known to be sharp \cite{13}.
In the case that $\OO$ is a $C^1$ domain, $\nabla (\CL)^{-1/2}$
is bounded on $L^p(\OO)$ for $1<p<\infty$.
We should point out that although our weighted
$L^2$ bounds are new, the boundedness of $\nabla (-\Delta)^{-1/2}$
on $L^p(\OO)$ for Lipschitz 
or $C^1$ domains was proved already in  \cite{13},
 by the method of complex
interpolation.
The direct extension of this method
to the case of continuous coefficients
fails, since it relies on the solvabilities of the $L^2$ 
Dirichlet and regularity problems.
However, Theorem B in this paper allows us to perturb
the operator $\CL$. This leads to the
following theorem.

\proclaim{\bf Theorem C} Let $\OO$ be a bounded Lipschitz domain
in $\br^n$, $n\ge 2$. Let $\CL$ be a second-order elliptic
operator of divergence form with real, symmetric, bounded
measurable coefficients on $\OO$, subject to Dirichlet
boundary condition.
Assume that the coefficients
$a_{jk}(x)$ are  in $VMO({\br^n})$.
Then there exists $\e>0$ such that
$\nabla (\CL)^{-1/2}$ is bounded on $L^p(\OO)$
for $1<p<3+\e$ if $n\ge 3$, and for $1<p<4+\e$ in the case
$n=2$. Consequently, there exists $\delta>0$ such that
$\nabla (\CL)^{-1/2}$ is bounded on 
 $L^2(\OO,\frac{dx}{\oo})$
where $\oo \in A_{\frac{4}{3}+\delta}(\br^n)$ if $n\ge 3$, and
 $\oo\in A_{\frac{3}{2}+\delta}(\br^2)$ in the case $n=2$.
If $\OO$ is a $C^1$ domain, $\nabla(\CL)^{-1/2}$
is bounded on $L^p(\OO)$ for any $1<p<\infty$, and
on $L^2(\OO, \oo dx)$ for any $\oo\in A_2(\br^n)$.
\endproclaim

\remark{\bf Remark 1.11}
For divergence form elliptic equations on $C^{1,1}$ domains
with $VMO$ coefficients,
the $W^{1,p}$ estimates were obtained in \cite{11} for any $1<p<\infty$.
The result was extended in \cite{3}
to the case of $C^1$ domains, for operators
with complex coefficients.
Our approach to Theorem C, which 
is very different from that in \cite{11,3},
is based on Theorem B and a perturbation argument found in \cite{6}.
Indeed, by Theorem B, it suffices to show that
solutions of $\CL u=0$ satisfies condition (i) in Theorem B
for some $p>3$ if $n\ge 3$, and $p>4$ in the case $n=2$.
To do this, we approximate $u$ on each ball by a solution
of a second order elliptic equation with constant
coefficients (Lemma 4.7). The desired estimate for $\nabla u$
follows from an approximation theorem (Theorem 4.13),
which is essentially proved in \cite{6}.
\endremark

The paper is organized as follows. Sections 2 and 3 are devoted
to the proof of Theorems A and B.
Theorem C is proved in section 4.
Finally in section 5 we give the proof of Theorems 3.1 and 3.2.

\noindent{\bf Acknowledgment.}
After this paper was submitted, the author was 
informed kindly by S.~Hofmann
of two recent preprints \cite{1,2} on the study of
Riesz transforms. In these two papers
 necessary and sufficient conditions are obtained
for the $L^p$
boundedness of Riesz transforms on manifolds \cite{2, Theorem 1.3}, and of
Riesz transforms associated to second order
elliptic operators with complex coefficients on $\br^n$
\cite{1, Proposition 5.6}.
The conditions are given in terms of the $L^p$ boundedness of
the operators $\sqrt{t}\nabla e^{-t\CL}$ uniformly for all $t>0$.
It is interesting to point out that
the key step in the proof of sufficiency of the conditions
in \cite{1,2}
uses a result similar to Theorem 3.1 of the present paper
(see Theorem 2.2 in \cite{1}).

The author thanks S.~Hofmann for pointing out the
relevance of the results in \cite{1,2}.
The author also would like to thank the referee for several 
valuable comments.

\bigskip

\centerline{\bf 2. Some Preliminaries}

In this section we will prove that condition (iii) in Theorems A
and B implies
(ii) which, in turn, implies condition (i).
We will also show that condition (i) leads to an $L^p$ estimate
on the kernel function of the resolvent $(\CL +\lambda)^{-1}$
for $\lambda
> 0$.

The following proposition is essentially due to Rubio de Francia \cite{16}.

\proclaim{\bf Proposition 2.1}
Let $T$ be a bounded operator on $L^2(E)$ where $E$
is a measurable subset of $\br^n$. Let $0<\delta\le 1$.
Suppose that
$$
\int_{E} |Tf|^2\,  \frac{dx}{\oo}
\le C\int_{E} |f|^2 \, \frac{dx}{\oo}
\ \ \text{ for any }\oo \in A_{1+\delta} (\br^n),
\tag 2.2
$$
where $C$ depends only on the $A_{1+\delta}$ bound of $\oo$.
Then $T$ is bound on $L^p(E)$ for $1<p<2/(1-\delta)$.
\endproclaim

\demo{Proof}
By considering the operator $\widetilde{T}(f)=\chi_E T(f\chi_E)$,
we may assume that $E=\br^n$.
It follows from assumption (2.2) that for any $\oo\in A_1$,
 $T$ is bounded on $L^2(\br^n,\frac{dx}{\oo})$ and 
 $L^2(\br^n, \oo^\delta dx)$.
It is known that the boundedness of $T$ on $L^2(\br^n, \frac{dx}{\oo})$
for any $\oo\in A_1$ implies its boundedness on 
$L^p(\br^n)$ for $1<p<2$, while the boundedness of $T$
on $L^2(\br^n, \oo^\delta dx)$ for any $\oo\in A_1$
implies its boundedness on $L^p(\br^n)$ for 
$2<p<2/(1-\delta)$.
We refer the reader to \cite{10, pp.141-142} for a simple and elegant
proof of these facts.
\enddemo

Using Proposition 2.1, it is easy to see that condition (iii)
in Theorem A or B implies condition (ii). 
Next we show that the $W^{1,p}$ estimate
follows from  the $L^p$ boundedness of the Riesz transform.

\proclaim{\bf Proposition 2.3}
Suppose that operator $\nabla (\CL)^{-1/2}$ is bounded
on $L^p(\br^n)$ for some $p>2$.
For $f\in L^p(\br^n)$ and $g\in L^q(\br^n)$ where $\frac{1}{p}
=\frac{1}{q}-\frac{1}{n}$, let $u\in W^{1,2}_{loc}(\br^n)$ 
be a weak solution of
$\CL u =\text{div} f +g$ in $\br^n$.
If $R^{\frac{n}{p}-\frac{n}{2}}
\left\{ \|\nabla u\|_{L^2(R\le |x|\le 2R)}
+R^{-1} \| u\|_{L^2(R\le |x|\le 2R)}\right\} \to 0$
as $R\to \infty$, then
$
\| \nabla u \|_p
\le C\big\{ \| f\|_p +\| g\|_q\big\}.
$
\endproclaim
\demo{Proof} The proof is rather standard. Let
$\varphi$ be a smooth cut-off function such that
$\varphi =1$ on $B(0,R)$, $\varphi =0$
outside of $B(0,2R)$, and $|\nabla \varphi |
\le C/R$. Then
$\CL(u \varphi) =\text{div}(\widetilde{f})
+\widetilde{g}$ in $\br^n$, where
$\widetilde{f}
=f\varphi -a_{jk} u \partial_k \varphi$ and
$\widetilde{g}=-f\nabla \varphi +g\varphi
-a_{jk}\partial_j u\partial_k \varphi$.
Since $u\varphi$, $\widetilde{f}$ and $\widetilde{g}$
all have compact supports, we may write
$u\varphi =\CL^{-1}( \text{div} \widetilde{f}
+\widetilde{g})$.

Suppose now that $\nabla(\CL)^{-1/2}$ is bounded on $L^p(\br^n)$
for some $p>2$. Since $\nabla (\CL)^{-1/2}$ is always
bounded on $L^t(\br^n)$ for any $1<t\le 2$,
it follows from duality 
that $\nabla \CL^{-1}\text{div} $ is bounded on $L^t(\br^n)$
for $2\le t\le p$.
This implies that
$$
\| \nabla (u\varphi)\|_t
 \le C\left\{\| \widetilde{f}\|_t
+\| (\CL)^{-1/2} \widetilde{g}\|_t\right\}
\le C 
\left\{\| \widetilde{f}\|_t
+\| \widetilde{g}\|_s\right\},
\tag 2.4
$$
where $\frac{1}{t}=\frac{1}{s}-\frac{1}{n}$.
We remark that the second inequality in (2.4) follows from the fact
that the kernel function $K(x,y)$ of the operator
$(\CL)^{-1/2}$ is bounded by $C|x-y|^{1-n}$.
Hence,
$$
\aligned
\| \nabla u\|_{L^t(B(0,R))}
&\le C\left\{ 
\| f\|_{L^t(B(0,2R))}
+\| g\|_{L^s(B(0,2R))}\right\}\\
&\ \ \ 
+C\, R^{-1}\left\{ \| u\|_{L^t(B(0,2R)\setminus B(0,R))}
+\| \nabla u\|_{L^s(B(0,2R)\setminus B(0,R))}
\right\}.
\endaligned
\tag 2.5
$$
By an iteration argument and Sobolev imbedding, this yields that
$$
\aligned
\| \nabla u\|_{L^p(B(0,R))}
&\le C\big\{ 
\| f\|_{L^p(B(0,CR))}
+\| g\|_{L^q(B(0,CR))}\big\}\\
&\ \ 
+C\, R^{\frac{n}{p}-\frac{n}{2}}
\big\{ R^{-1}\, \| u\|_{L^2(R\le |x|\le CR)}
+\| \nabla u\|_{L^2(R\le |x|\le CR)}
\big\}.
\endaligned
\tag 2.6
$$
Letting $R\to \infty$ in (2.6), one obtains the desired estimate.
\enddemo

\proclaim{\bf Lemma 2.7}
In Theorem A or B, condition (ii) implies condition (i).
\endproclaim

\demo{Proof}
Let $p>2$, and suppose that $\nabla (\CL)^{-1/2}$ is bounded on $L^p(\br^n)$.
Then the operator is bounded on $L^t(\br^n)$ for $1<t\le p$.
Let $u$ be a weak solution of $\CL u=0$
in $3B=B(x_0,3r)$. For $1<\gamma_1<\gamma_2<3/2$, 
choose a smooth cut-off function $\varphi$ such that
$\varphi=1$ on $\gamma_1 B$, $\varphi =0$ outside of
$\gamma_2 B$, and $|\nabla\varphi|\le C/r$.
Note that $\CL(u\varphi)
=-\partial_j\left( a_{jk} u \partial_k \varphi\right)
-a_{jk} \partial_k u \partial_j \varphi$
in $\br^n$. By Proposition 2.3, if $|\nabla u|\in L^s(\gamma_2B)$,
then $|\nabla u|\in L^t(\gamma_1 B)$ and
$$
\|\nabla u\|_{L^t(\gamma_1 B)}
\le C \, r^{-1} \left\{  \| u\|_{L^t(\gamma_2 B)}
+\| \nabla u\|_{L^s(\gamma_2 B)}\right\},
\tag 2.8
$$
where $2<t\le p$ and $\frac{1}{t}=\frac{1}{s}-\frac{1}{n}$.
Since $u-c$ is also a weak solution in $3B$, we may use the
 $L^\infty$  estimate and Poincar\'e inequality to obtain
$$
\aligned
\bigg(\frac{1}{|B|}
&\int_{\gamma_1 B} |\nabla u|^t dx\bigg)^{1/t}\\
&\le  C\,
\left\{ 
\left(\frac{1}{|2B|}
\int_{2B} |\nabla u|^2 dx\right)^{1/2}
+\left(\frac{1}{|B|}
\int_{\gamma_2 B} |\nabla u|^s dx\right)^{1/s}
\right\}.
\endaligned
\tag 2.9
$$
From this, estimate (1.3) in Theorem A
follows by an iteration argument, starting with $s=2$.

In the case of Theorem B, we first choose $r_0>0$ so small that
for any $P\in \partial\OO$, $\OO\cap B(P,3r_0)$
is given by the intersection of 
the region above a Lipschitz graph and $B(P,3r_0)$,
 after a possible rotation of the coordinate system.
Given $B(x_0,r)$ with $0<r<r_0$, consider two cases:
(1) $B(x_0,3r)\subset \OO$,
(2)  $x_0\in\partial\OO$.
The first case may be treated exactly as in Theorem A.
In the second case, instead of replacing $u$ by $u-c$ and
using Poincar\'e inequality, one applies the Poincar\'e
inequality on $\OO\cap B(x_0,\gamma_2 r)$
for functions which vanish on $B(x_0, 3r)\cap \partial\OO$.
The rest is the same.
\enddemo

To complete the proof of Theorems A and B, it remains to show that
condition (i) implies conditions (ii) and (iii). To this end, 
we need to
estimate the kernel function $\Gamma_\lambda (x,y)$ of the
resolvent $(\CL +\lambda)^{-1}$ for $\lambda> 0$. We begin
with a size estimate for $n\ge 3$:
$$
|\Gamma_\lambda(x,y)|
\le C\, e^{-c \sqrt{\lambda}|x-y|}\cdot
\frac{1}{|x-y|^{n-2}},
\tag 2.10
$$
which follows directly from the formula
$(\CL +\lambda)^{-1}
=\int_0^\infty e^{-\lambda t} e^{-t\CL}\, dt$ and
 the well known upper bound for the heat kernel of $\CL$ \cite{9}.
In the case $n=2$, one needs to replace $\frac{1}{|x-y|^{n-2}}$
by $|\ln (\sqrt{\lambda}|x-y|)| +1$.
The rest of this section is devoted to the proof of the
following theorem. We remark that estimates similar to
(2.12)-(2.13) may be found in \cite{10}.

\proclaim{\bf Theorem 2.11}
Suppose that operator $\CL$ in (1.1) satisfies condition
(i) in Theorem B for some $p>2$. Then, if $n\ge 3$,
$$
\align
& \left(\frac{1}{r^n}
\int_{\{ x\in \OO:\ r\le |x-y|\le 2r\} }
|\nabla_x \Gamma_\lambda (x,y)|^p
dx\right)^{1/p}
\le C\, e^{-c \sqrt{\lambda} r}\cdot\frac{1}{r^{n-1}},
\tag 2.12\\
&
\left(\frac{1}{r^n}
\int_{\{ x\in\OO:\ r\le |x-y|\le 2r\} }
|\nabla_x \Gamma_\lambda (x,y)-\nabla_x \Gamma_\lambda
(x,y+h)|^p
dx\right)^{1/p}\\
& \ \ \ \ \ \ \
\le C\, \left(\frac{|h|}{r}\right)^\eta
\cdot
 e^{-c \sqrt{\lambda} r}
\cdot\frac{1}{r^{n-1}},
\tag 2.13
\endalign
$$
where $0<r<c\,r_0$,
 $y, y+h\in \OO$, $|h|<c\,r$, and $\eta
=\eta (n,\mu,\OO)>0$. If $n=2$, one needs to replace
$\frac{1}{r^{n-1}}$ in (2.12)-(2.13)
by $\big(|\ln (\sqrt{\lambda} r)| +1\big)/r$.
If $\OO=\br^n$ and $\CL$ satisfies condition (i) in Theorem A, 
above estimates hold for any $0<r<\infty$.
\endproclaim

Note that $\Gamma_\lambda (x,y)
=\Gamma_\lambda(y,x)$, and
$\Gamma_\lambda
(\cdot,y)$ is a weak solution of $\CL u +\lambda u =0$
in $\OO\setminus \{ y\}$. Using size estimate (2.10) and
H\"older estimates, it is easy to see that
Theorem 2.11 is a consequence of the following lemma.

\proclaim{Lemma 2.14}
Assume that $\CL$ satisfies condition (i) in Theorem B
for some $p>2$. Then there exist constants $r_1>0$,
 $\alpha>1$ and $C>0$ independent of $\lambda>0$,
 such that
if $u$ is a weak solution of $\CL u+\lambda u=0$ in $
B(x_0, \alpha r)\cap\OO$ for some $x_0\in\overline{\OO}$,
$0<r<r_1$ and $u=0$ on $B(x_0,\alpha r)\cap\partial\Omega$, then
$$
\left(\frac{1}{r^n}\int_{B(x_0,r)\cap\OO} |\nabla u|^p dx\right)^{1/p}
\le \frac{C}{r}
\left(\frac{1}{r^n}\int_{B(x_0,\alpha r)\cap\OO}
| u|^2 dx\right)^{1/2}
\tag 2.15
$$
If $\OO=\br^n$ and $\CL$ satisfies condition (i) in Theorem A, 
above statement holds for $r_1=\infty$.
\endproclaim

\demo{Proof} Let $u$ be a weak solution of $\CL u +\lambda u=0$
in $B(x_0,\alpha r)\cap\OO$ and
$u=0$ on $B(x_0,\alpha r)\cap\partial\Omega$,
 where $\alpha =2\alpha_2$.
We only consider the case $x_0\in \partial\OO$.

Let $D=B(x_0,r)\cap \OO$ and $tD=B(x_0,tr)\cap\OO$.
Let $v$ be a weak solution of $\CL v =0$ in $\alpha_2 D$ such that
$w\equiv u-v \in H_0^1(\alpha_2 D)$. Using condition (i), we have
$$
\left(\frac{1}{r^n}
\int_D |\nabla u|^p dx\right)^{1/p}
\le C \left( \frac{1}{r^n}\int_{\alpha_1 D} |\nabla u|^2 dx\right)^{1/2}
+C \left(\frac{1}{r^n}\int_{\alpha_1 D}
|\nabla w|^p dx\right)^{1/p}.
\tag 2.16
$$
To estimate $\nabla w$ on $\alpha_1 D$, observe 
 that $\CL w=-\lambda u$ in $\alpha_2 D$. Hence we may write
$$
w(x)=-\lambda \int_{\alpha_2 D} G(x,y) u(y) dy,
\tag 2.17
$$
where $G(x,y)$ is the Green's function for
$\CL$ on $\alpha_2 D$. It follows that
$$
\aligned
|\nabla w (x)| &\le \lambda \| u\|_{L^\infty (\alpha_2 D)}
\int_{\alpha_2 D} |\nabla_x G(x,y)| dy\\
&
\le \frac{C}{r^2}\left(\frac{1}{r^n}
\int_{2\alpha_2 D} |u|^2 dx \right)^{1/2} h(x),
\endaligned
\tag 2.18
$$ 
where $h(x)=\int_{\alpha_2 D} |\nabla_x G(x,y)| dy$, and we have used the
Cacciopoli inequality
$$
\lambda \, \int_{\frac32 \alpha_2 D}
|u|^2\, dx
+\int_{\frac32\alpha_2 D} |\nabla u|^2 \, dx
\le
\frac{C}{r^2}
\int_{2\alpha_2 D} |u|^2\, dx.
\tag 2.19
$$ Note that
$$
\| h\|_{L^p(\alpha_1 D)}
=\sup_{\| g\|_{p^\prime}\le 1}
\bigg| \int_{\alpha_1 D} h(x) g(x) dx\bigg|
\le \sup_{\| g\|_{p^\prime}\le 1}
\int_{\alpha_2 D} |Tg(y)| dy,
\tag 2.20
$$
where 
$$
Tg(y)=\int_{\alpha_1 D} |\nabla_x G(x,y)| g(x) dx.
\tag 2.21
$$
Since $G(\cdot, y)$ is a weak solution of
$Lu=0$ in $\alpha_2 D\setminus\{ y\}$, it follows
from estimate (1.10) that
$$
\left(\int_{E_j} |\nabla_x G(x,y)|^p dx \right)^{1/p}
\le C\, (2^{-j} r)^{1-\frac{n}{p^\prime}},
\tag 2.22
$$
where $E_j=E_j(y)=\left\{ x\in \alpha_1 D: 2^{-j-1} r\le |x-y|\le 2^{-j} r
\right\}$ for $j\ge -3$. Thus, by H\"older inequality,
$$
\aligned
|Tg(y)|
&\le \sum_{j=-3}^\infty
\left(\int_{E_j}
|\nabla_x G(x,y)|^p dx\right)^{1/p}
\left(\int_{E_j} |g|^{p^\prime} dx\right)^{1/p^\prime}\\
&\le C\,r \left\{ M(|g|^{p^\prime})(y)\right\}^{1/p^\prime},
\endaligned
\tag 2.23
$$
where $M$ denotes the Hardy-Littlewood maximal function.
By Kolmogorov's Lemma \cite{10, p.102}, this implies that
$$
\int_{\alpha_2 D} |Tg(y)| dy
\le C\, r |\alpha_2 D|^{\frac{1}{p}} \| g\|_{p^\prime}.
\tag 2.24
$$
In view of (2.18) and (2.20), we obtain
$$
\aligned
\|\nabla w\|_{L^p(\alpha_1 D)}
&\le \frac{C}{r^2} \left(\frac{1}{r^n}\int_{2\alpha_2 D}
|u|^2 dx\right)^{1/2} \| h\|_{L^p (\alpha_1 D)}\\
&\le C\, r^{\frac{n}{p}-1}
\left(\frac{1}{r^n}\int_{2\alpha_2 D} |u|^2 dx\right)^{1/2}.
\endaligned
\tag 2.25
$$
The desired estimate (2.15) with $\alpha
=2\alpha_2$ now follows from (2.16) and (2.25).
\enddemo

\bigskip

\centerline{\bf 3. Proof of Theorems A and B}

In this section we show that condition (i) in Theorem A or B
implies conditions (ii) and (iii). This,
together with Proposition 2.1 and Lemma 2.7,
 completes the proof of Theorems A
and B.

The proof of condition (i) implying (ii) relies on 
Theorem 3.1, which may be considered as a refined (and dual) version of
the well known Calder\'on-Zygmund Lemma. 
Its proof as well as the proof of Theorem 3.3 will be
given in section 5.

\proclaim{\bf Theorem 3.1}
Let $T$ be a bounded sublinear operator on $L^2(\br^n)$. 
Let $p>2$.
Suppose that there exist constants $\alpha_2>\alpha_1>1$, $N>1$ such that
$$
\aligned
&\left\{ \frac{1}{|B|}
\int_B |Tf|^p dx \right\}^{1/p}\\
&
\le 
N\, \left\{ \left(\frac{1}{|\alpha_1 B|}
\int_{\alpha_1 B} |Tf|^2 dx\right)^{1/2}
+ \sup_{B^\prime \supset B}
\left(\frac{1}{|B^\prime|}
\int_{B^\prime}
|f|^2 dx\right)^{1/2}\right\},
\endaligned
\tag 3.2
$$
for any ball $B
\subset \br^n$, and any bounded measurable function $f$
 with compact 
supp$(f)\subset \br^n\setminus \alpha_2 B$.
Then $T$ is bounded on $L^q(\br^n)$ for any $2<q<p$.
\endproclaim

Theorem 3.1 may be extended to the case of bounded Lipschitz domains.

\proclaim{\bf Theorem 3.3}
Let $T$ be a bounded sublinear operator on $L^2(\OO)$, where $\OO$ is 
a bounded Lipschitz domain in $\br^n$. Let $p>2$. Suppose that there exist
constants $r_0>0$, $N>1$ and $\alpha_2>\alpha_1>1$ such that for any bounded
measurable function $f$ with
supp$(f)\subset \OO\setminus \alpha_2 B$,
$$
\aligned
&\left\{
\frac{1}{r^n}\int_{\OO\cap B}
|Tf|^p\, dx\right\}^{1/p}\\
&\le N
\left\{
\left(\frac{1}{r^n}\int_{\OO\cap \alpha_1 B}
|Tf|^2\, dx\right)^{1/2}
+\sup\Sb B^\prime\supset B\endSb
\left(\frac{1}{|B^\prime|}
\int_{B^\prime} |f|^p\, dx\right)^{1/p}\right\},
\endaligned
\tag 3.4
$$
where $B=B(x_0,r)$ is a ball with the property that
$0<r<r_0$ and either $x_0\in \partial\OO$ or
$B(x_0,\alpha_2 r)\subset \OO$.
Then $T$ is
bounded on $L^q(\OO)$ for any $2<q<p$.
\endproclaim

\proclaim{Lemma 3.5} In Theorem A or B, condition (i)
implies (ii).
\endproclaim

\demo{Proof} We first consider the case $\OO=\br^n$.
Assume that operator $\CL$ satisfies condition (i) in Theorem A
for some $p>2$. By Theorem 3.1, the linear operator
$T=\nabla (\CL)^{-1}\text{div}$ is bounded on
$L^q(\br^n)$ for $2\le q<p$. Indeed,
$T$ is clearly bounded on $L^2(\br^n)$.
To verify (3.2), we let $u=(\CL)^{-1}\text{div} (f)$ where $f$ is a
bounded measurable function
with compact supp$(f)\subset \br^n\setminus 4B$. Observe that
$\CL u =0$ in $3B$. Thus inequality (3.2) follows directly from
 condition (i). By Theorem 3.1 and duality,
$\nabla (\CL)^{-1}\text{div} $ is bounded on $L^q$
for $p^\prime<q<p$.

Next, since $\|L^{1/2} f\|_q\le C\| \nabla f\|_q $ for any $1<q<\infty$
\cite{4, p.114}, we 
have
$$
\| (\CL)^{-1/2}\text{div} f\|_q =\| \CL^{1/2} (\CL)^{-1} \text{div} f \|_q
\le C\, \| \nabla (\CL)^{-1} \text{div} f\|_q,
\le C \, \| f\|_q
\tag 3.6
$$
where $p^\prime<q<p$. Consequently, by duality,
 $\nabla (\CL)^{-1/2}$
is bounded on $L^q(\br^n)$ for $1<q<p$.
Finally by the self-improvement property of condition (i)
(see Remark 1.9),
 we may conclude that $\nabla (\CL)^{-1/2}$ is bounded
on $L^q(\br^n)$ for $1<q<p+\e$.

The proof is similar in the case of Theorem B. In the place of
Theorem 3.1, we use Theorem 3.3. Also we note that 
for a bounded Lipschitz domain,
 the inequality
$\| \CL^{1/2} f\|_q\le C\, \| \nabla f\|_q$
has been established in \cite{5}.
 The proof is finished.
\enddemo

To show that condition (i) implies the $L^2$ weighted norm inequality
for the Riesz transform, we use the functional calculus formula
$$
(\CL)^{-1/2}=\frac{1}{\pi}
\int_0^\infty \lambda^{-1/2} (\CL +\lambda)^{-1} d\lambda
\tag 3.7
$$ 
to write
$$
\nabla (\CL)^{-1/2} f(x) =\int_{\OO} K(x,y) f(y)\, dy,
\tag 3.8
$$
where
$$
K(x,y)=\frac{1}{\pi}\int_0^\infty \lambda^{-1/2} \nabla_x \Gamma_\lambda
(x,y) d\lambda,
\tag 3.9
$$
and $\Gamma_\lambda (x,y)$ is the Green's function
 for $\CL +\lambda$.

\proclaim{Lemma 3.10}
Suppose that $\CL$ satisfies condition (i) in Theorem B
for some $p>2$. Then
$$
\align
\left(\frac{1}{r^n}
\int_{\{ x\in\OO:\, r\le |x-y|\le 2r\} }
|K(x,y)|^p \, dx \right)^{1/p}
& \le \frac{C}{r^n},
\tag 3.11\\
\left(\frac{1}{r^n}
\int_{\{ x\in\OO:\, r\le |x-y|\le 2r\} }
|K(x,y)-K(x,y+h)|^p \, dx \right)^{1/p}
&\le
\frac{C}{r^n}\cdot \left(\frac{|h|}{r}\right)^\eta
\tag 3.12
\endalign
$$
where $0<r<r_1$,
$y,y+h\in \OO$, $|h|\le c\,r$, and $\eta=\eta(n,\mu,\OO)>0$.
If $\OO=\br^n$ and $\CL$ satisfies condition (i) in Theorem A,
then estimates (3.11)-(3.12) hold for all $0<r<\infty$.
\endproclaim

\demo{Proof} In view of  (3.9), 
estimates
(3.11) and (3.12) follow directly from (2.14) and (2.15) respectively
by integration.
\enddemo

To use the estimates in Lemma 3.10 effectively, we consider the 
adjoint operator of $\nabla (\CL)^{-1/2}$:
$$
Sf(x)=\int_{\OO}
K(y,x) f(y)\, dy.
\tag 3.13
$$
Recall that the sharp function of $f$ is defined by
$$
f^\# (x)
\equiv
\sup_{B\owns x}\inf_{\beta\in\br}
\frac{1}{|B|}
\int_{B}| f(y) -\beta| dy.
\tag 3.14
$$

\proclaim{Lemma 3.15}
If operator $\CL$ satisfies condition (i) in Theorem A for some
$p>2$, then
$$
(Sf)^\# (x)\le C\, \left\{ M(|f|^{p^\prime})(x)\right\}^{1/p^\prime}
\tag 3.16
$$
for any $x\in \br^n$. If $\CL$ satisfies condition (i) in
Theorem B, and $Sf$ is defined to be zero outside of $\OO$,
then (3.16) holds for any $x\in\OO$.
\endproclaim

\demo{Proof} We first consider the case of Theorem A.
Suppose $x\in B=B(x_0,r)$. Let $f=g+h$ where $g=f \chi_{4B}$. 
Since $S$ is bounded on $L^{p^\prime}(\br^n)$ by Lemma 3.4 and duality, 
we have
$$
\aligned
\frac{1}{|B|}\int_B |S (g)| dy
&
\le \left(\frac{1}{|B|}\int_B |S(g)|^{p^\prime} dy\right)^{1/p^\prime}\\
&\le C\left(\frac{1}{|B|}
\int_{4B} |f|^{p^\prime} dy\right)^{1/p^\prime}
\le C\left\{ M(|f|^{p^\prime})(x)\right\}^{1/p^\prime}.
\endaligned
\tag 3.17
$$

Next, let $\beta =S(h)(x_0)$. It follows from H\"older inequality
and estimate (3.12) that, for $y\in B$,
$$
\aligned
&
|S(h)(y)-\beta|\le
\int_{\br^n\setminus 4B}
|K(z,y)-K(z,x_0)|\, |f(z)| dz\\
&
\le \sum_{j=2}^\infty
\left(
\int_{2^{j+1}B\setminus 2^j B}
|K(z,y)-K(z,x_0)|^{p} dz\right)^{1/p}
\left(
\int_{2^{j+1}B\setminus 2^j B}
|f(z)|^{p^\prime} dz \right)^{1/p^\prime} \\
&\le C\, \left\{ M (|f|^{p^\prime})(x)\right\}^{1/p^\prime}
\sum_{j=2}^\infty
2^{-j\eta}
\le C \left\{ M (|f|^{p^\prime})(x)\right\}^{1/p^\prime}.
\endaligned
$$
This, together with (3.17), gives
$$
\frac{1}{|B|}\int_B
|Sf (y)-\beta| dy
\le 
C \left\{ M (|f|^{p^\prime})(x)\right\}^{1/p^\prime},
\tag 3.18
$$
from which (3.16) follows.

In the case of Theorem B, we may use the same argument as above to
show that for any $x\in \OO$, estimate (3.18) holds for
any ball $B=B(x_0,r)\owns x$
 with $0<r<r_1$. If $r\ge r_1$, we use the boundedness
of $S$ on $L^{p^\prime}(\OO)$ to obtain
$$
\aligned
\frac{1}{|B|}\int_{B\cap\OO} |Sf|\, dy
&\le C_{r_1}\, \left\{\int_{\OO} |Sf|^{p^\prime}\, dy\right\}^{1/p^\prime}
\le C\left\{ \int_{\OO}
|f|^{p^\prime}\, dy\right\}^{1/p^\prime}\\
&\le C\left\{ M(|f|^{p^\prime})(x)\right\}^{1/p^\prime}.
\endaligned
\tag 3.19
$$
The proof is complete.
\enddemo

\proclaim{\bf Proposition 3.20}
Let $f\in L^1_{loc} (E)$ where $E$ is a measurable set of $\br^n$.
Suppose $\oo\in A_2(\br^n)$ and $f\in L^2(E, \oo dx)$. Define
$f$ to be zero on $\br^n\setminus E$. Then
$$
\int_E |f(x)|^2 \, \oo \, dx
\le C\, \int_E |f^\# (x)|^2 \, \oo \, dx.
\tag 3.21
$$
where $C$ depends only on $n$ and the $A_2$ bound of $\oo$.
\endproclaim

\demo{Proof} The case $E=\br^n$ is well known. It was proved in 
\cite{10}, using the following good-$\lambda$ inequality
$$
\oo\big\{ x\in \br^n:\,
M_d f(x)>2\lambda,\ f^\#(x)\le \gamma\lambda\big\}
\le C\, \gamma^\delta \oo\big\{
x\in \br^n: \, M_df(x)>\lambda\big\}
\tag 3.22
$$
where $\oo\in A_\infty$, and
$M_df$ denotes the dyadic maximal function of $f$.
In general, we use (3.22) to obtain
$$
\aligned
&\oo\big\{ x\in E: \, M_df(x)>\lambda\big\}\\
&\le \oo\big\{ x\in E: \, f^\#(x)>\gamma \lambda\big\}
+C\, \gamma^\delta \oo
\big\{ x\in \br^n:\ M_df(x)>\lambda\big\}.
\endaligned
\tag 3.23
$$
By integration, this gives
$$
\aligned
\int_E |f(x)|^2\, \oo\, dx
& \le C_\gamma\, \int_E |f^\# (x)|^2 \, \oo\, dx
+C\, \gamma^\delta \int_{\br^n}
|M_d f(x)|^2\, \oo\, dx\\
&
\le C_\gamma\, \int_E |f^\# (x)|^2\, \oo\, dx
+C\, \gamma^\delta\, \int_E |f(x)|^2\, \oo\, dx,
\endaligned
\tag 3.24
$$
where we have used $\oo\in A_2$ and the weighted norm inequality for 
$M_d$. Inequality (3.21) now follows by choosing
$\gamma$ so small that $C\gamma^\delta <1/2$.
\enddemo

We are now in a position to complete the proof of Theorems A and B.

\proclaim{Lemma 3.25}
In Theorem A or B, condition (i) implies (iii).
\endproclaim

\demo{Proof} We give the proof for the case of Theorem A.
The case of Theorem B is similar.

Suppose that operator $\CL$ satisfies condition (i) in Theorem A
for some $p>2$.
 It follows from
 Proposition 3.20, Lemma 3.15 and  the weighted norm inequality for
the Hardy-Littlewood maximal function that
$$
\aligned
\int_{\br^n} |Sf|^2 \,\oo dx
&\le C\, \int_{\br^n} |(Sf)^\# |^2 \, \oo dx
\le C\, \int_{\br^n}
\left\{ M(|f|^{p^\prime})\right\}^{2/p^\prime} \, \oo dx\\
&
\le C\, \int_{\br^n} |f|^2 \, \oo dx,
\endaligned
\tag 3.26
$$
where $f$ is a bounded function with compact support, and
$\oo\in A_{2/p^\prime} (\br^n)$.
We remark that the first inequality in (3.26) requires
$Sf\in L^2(\br^n, \oo dx)$. To see this, let us assume that
supp$(f)\subset B(0,R)$. Since $S$ is bounded on 
$L^q(\br^n, dx)$ for $q\ge p^\prime$, 
$Sf\in L^2(B(0,2R),\oo dx)$ by H\"older inequality.
If $|x|\ge 2R$, using (3.11), we may show that
$$
|Sf(x)|\le C\, \| f\|_{p^\prime} \cdot \frac{1}{|x|^{n(1-\frac{1}{p})}}.
\tag 3.27
$$
This is enough to assure that $Sf\in L^2(\br^n\setminus B(0,2R), \oo dx)$
for any $\oo\in A_{2/p^\prime}$.

Finally, by (3.26) and
 duality, $\nabla (\CL)^{-1/2}$ is bounded on $L^2(\br^n, \frac{dx}
{\oo})$ for any $\oo \in A_{2/p^\prime}(\br^n)$.
Condition (iii) now follows by the self improvement
property of condition (i).
\enddemo

\bigskip

\centerline{\bf 4. Operators with VMO Coefficients on Lipschitz Domains}
 
In this section we will prove Theorem C stated in the
Introduction. To do this, we first show that operators
with constant coefficients satisfy condition (i)
in Theorem B. We then use an  approximation argument
found in \cite{6} to prove that operators with VMO
coefficients also satisfy condition (i).

\proclaim{\bf Lemma 4.1}
Suppose operator $\CL$ in (1.1) has constant coefficients.
Then
 it satisfies condition (i) in Theorem B for some $p>3$ if $n\ge 3$,
and for some $p>4$ in the case $n=2$. Moreover,
$p$ depends only on $\OO$, $n$ and the ellipticity constant
$\mu$ of $\CL$.
\endproclaim

\demo{Proof} We may assume that $\CL=-\Delta$.
If $B(x_0,3r)\subset \OO$, inequality (1.5) for any $p>2$
follows easily from
the interior estimates. Suppose $x_0\in \partial\OO$.
We may assume that $x_0=0$ and 
$$
D(0,r_1)=\OO\cap B(0,r_1)
=\big\{ (x^\prime,x_n)\in \br^n:\
x_n>\psi(x^\prime)\big\}
\cap B(0,r_1),
\tag 4.2
$$
where $\psi(x^\prime)$ is a Lipschitz function on $\br^{n-1}$.

Let $\alpha_2>\alpha_1>1$.
Suppose $0<r<r_0=c\, r_1$, where $c>0$ is sufficiently small.
Let $u\in H^1(D(0,\alpha_2 r))$ be a harmonic function
in $D(0,\alpha_2  r)$ such that $u=0$ on $B(0,\alpha_2 r)\cap\partial\OO$.
By the boundary H\"older estimate and Poincar\'e inequality, 
for $x=(x^\prime,x_n)
\in D(0,r)$, we have
$$
|\nabla u (x)|
\le 
C\, \big(x_n-\psi(x^\prime)\big)^{\eta -1}\cdot\frac{1}{r^{\eta -1}}\cdot
\left\{\frac{1}{r^n}
\int_{D(0, \alpha_1 r)} |\nabla u|^2\, dy\right\}^{1/2},
\tag 4.3
$$
where $\eta >1/2$ if $n=2$, and $\eta>0$ if $n\ge 3$.
It follows that
if $(p-2)(\eta -1)>-1$,
$$
\aligned
&\int_{D(0,r)}
|\nabla u|^p\, dx\\
&\le \frac{C}{r^{(p-2)(\eta -1)}}
\int_{D(0,r)}
|\nabla u|^2 \, \big(x_n-\psi(x^\prime)\big)^{(p-2)(\eta -1)} \, dx
\cdot\left\{\frac{1}{r^n}\int_{D(0,\alpha_1 r)}
|\nabla u|^2 dx\right\}^{\frac{p}{2}-1}\\
&
\le \frac{C}{r^{(p-2)(\eta -1)}}
\int_{\Delta_r}
|(\nabla u)^*|^2 \, d\sigma
\cdot
\int_0^{cr}
 t^{(p-2)(\eta -1)}\, dt
\cdot\left\{\frac{1}{r^n}\int_{D(0,\alpha_1 r)}
|\nabla u|^2 dy\right\}^{\frac{p}{2}-1}\\
&
\le C\, r\, \int_{\Delta_r} |(\nabla u)^*|^2\, d\sigma
\cdot
\left\{\frac{1}{r^n}\int_{D(0,\alpha_1 r)}
|\nabla u|^2 dx \right\}^{\frac{p}{2}-1},
\endaligned
$$
where $\Delta_r=\{ (x^\prime,\psi(x^\prime)):\, |x^\prime|<r\}$
and $(\nabla u)^* (x^\prime, \psi(x^\prime))
=\sup\{ |\nabla u(x^\prime, x_n)|:\ (x^\prime, x_n)\in D(0,r)\}$.
This, together with the inequality
$$
\int_{\Delta_r} |(\nabla u)^*|^2 \, d\sigma
\le \frac{C}{r}
\int_{D(0,\alpha_1 r)}
|\nabla u|^2 \, dx,
\tag 4.4
$$
gives the desired estimate (1.5) for $2<p<\bar{p}$, where
$\bar{p}=2+\frac{1}{1-\eta}$. Note that $\bar{p}>3$ if
$n\ge 3$, and $\bar{p}>4$ if $n=2$.
Finally we point out that since $u=0$ on $\Delta_{\alpha_2 r}$,
(4.4) follows from the $L^2$ solvability of the 
regularity problem for Laplace's equation on Lipschitz domains,
by an integration argument (see \cite{8}).
\enddemo

\remark{\bf Remark 4.5}
If $\OO$ is a $C^1$ domain, then H\"older estimate (4.3)
holds for any $0<\eta <1$. It follows that operators with constant
coefficients satisfy condition (i) in Theorem B 
for any $p>2$. Consequently, the Riesz transform
$\nabla (\CL)^{-1/2}$ is bounded on $L^p(\OO)$ for all
$1<p<\infty$.
\endremark

A function $f$ in $BMO(\br^n)$ is said to be in $VMO(\br^n)$
if
$$
\lim_{r\to 0}\sup_{x_0\in \br^n}
\frac{1}{r^n}
\int_{B(x_0,r)}
|f- f_{B(x_0,r)}| \, dx
=0,
\tag 4.6
$$
where $f_{B(x_0,r)}=\int_{B(x_0,r)} f dx/|B(x_0,r)|$ is the average
of $f$ over $B(x_0,r)$.

\proclaim{\bf Lemma 4.7}
Let $\OO$ be a bounded Lipschitz domain in $\br^n$. Suppose
that the coefficients of operator $\CL$ in (1.1) are in 
$VMO(\br^n)$. Then there exist a function
$\phi(r)$ and some constants $C>0, \alpha>1$, $r_1>0$ and
$p>3$ ($p>4$, if $n=2$)
with the following properties:
(1) $\lim_{r\to 0}\phi(r)=0$, (2)
for any weak solution of $\CL u=0$ in $D(x_0,\alpha r)$ and
$u=0$ on $B(x_0,\alpha r)\cap \partial\OO$ with
$x_0\in \overline{\OO}$ and $0<r<r_1$, there exists
a function $v\in W^{1,p}(D(x_0,r))$ such that
$$
\align
\left\{ \frac{1}{r^n}
\int_{D(x_0,r)}
|\nabla v|^p \, dx\right\}^{1/p}
& \le C\, \left\{ \frac{1}{r^n}
\int_{D(x_0,\alpha r)}
|\nabla u|^2 dx\right\}^{1/p},\tag 4.8\\
\left\{ \frac{1}{r^n}
\int_{D(x_0,r)}
|\nabla u -\nabla v|^2\, dx\right\}^{1/2}
&\le \phi(r)
\left\{\frac{1}{r^n}
\int_{D(x_0,\alpha r)}
|\nabla u|^2\, dx\right\}^{1/2}.\tag 4.9
\endalign
$$
\endproclaim

\demo{Proof}
Let $u$ be a weak solution of $\CL u=0$ in $D(x_0, \alpha r)$
and $u=0$ on $B(x_0,\alpha r)\cap \partial\OO$ with
$x_0\in \overline{\OO}$ and $0<r<r_1$.
Consider the operator $\CL_0=-\partial_j b_{jk}\partial_k$, 
where $b_{jk}$ is a
constant given by
$$
b_{jk}=\frac{1}{|B(x_0, \alpha r)|}
\int_{B(x_0,\alpha r)} a_{jk} (x)\, dx.
\tag 4.10
$$
Let $v$ be a weak solution of $\CL_0 v=0$ in $D(x_0,\beta r)$
such that $u-v\in H^1_0(D(x_0,\beta r))$, where $\beta=\alpha/2=\alpha_1$.
We will show that $v$ satisfies estimates (4.8)-(4.9).
To this end, we first note that  $v=0$ on $B(x_0,\beta r)\cap
\partial\OO$. Thus, by Lemma 4.1,
$$
\aligned
\left\{\frac{1}{r^n}
\int_{D(x_0,r)}
|\nabla v|^p \, dx\right\}^{1/p}
&\le C\left\{\frac{1}{r^n}
\int_{D(x_0, \beta r)} |\nabla v|^2\, dx\right\}^{1/2}\\
&\le C\, \left\{\frac{1}{r^n}
\int_{D(x_0, \beta r)} |\nabla u|^2\, dx\right\}^{1/2},
\endaligned
\tag 4.11
$$
where $p>3$ for $n\ge 3$, and $p>4$ if $n=2$.
This gives (4.8).

To see (4.9), we observe that
$\CL_0(u-v)
=(\CL_0 -\CL) u =-\partial_j (b_{jk}-a_{jk})\partial_k u$.
It follows from the energy estimate that
$$
\aligned
&\left\{\frac{1}{r^n}
\int_{D(x_0, r)}
|\nabla u-\nabla v|^2\, dx\right\}^{1/2}\\
&  \le C\, \sum_{j,k}
\left\{ \frac{1}{r^n}\int_{D(x_0, \beta r)}
|b_{jk}-a_{jk}|^2\, |\nabla u|^2\, \, dx\right\}^{1/2}\\
&
\le C\, \sum_{j,k}
\left\{\frac{1}{r^n}\int_{B(x_0,\beta r)} |b_{jk}-a_{jk}|^{2q^\prime} dx
\right\}^{1/(2q^\prime)}
\left\{\frac{1}{r^n}
\int_{D(x_0,\beta r)}
|\nabla u|^{2q}dx\right\}^{1/(2q)}\\
&
\le \phi(r)\,
\left\{\frac{1}{r^n}
\int_{D(x_0,\alpha r)}
|\nabla u|^{2}dx\right\}^{1/2}
\endaligned
$$
where $q=1+\delta$ and $\delta>0$ is so small that
the $L^{2q}$ estimates hold for solutions of $\CL u=0$
\cite{15}.
Also we have introduced the function  $\phi (r)$ by
$$
\phi(r)=
C\, \sup_{x_0\in\overline{\OO}}
\sum_{j,k}
\left\{\frac{1}{r^n}\int_{B(x_0,\alpha r)}
|a_{jk}-b_{jk}|^{2q^\prime}\, dx\right\}^{1/(2q^\prime)}.
\tag 4.12
$$
Finally we note that by the John-Nirenberg inequality, if
$a_{jk}\in VMO(\br^n)$, then
$\phi(r)\to 0$ as $r\to 0$.
This completes the proof.
\enddemo

With Lemma 4.7 at our disposal, we may invoke the following
approximation theorem to finish the proof of
Theorem C. 

\proclaim{\bf Theorem 4.13}
Let $f:E\to\br^m$ be a locally square integrable function,
where $E$ is an open set of $\br^n$. Let $p>2$.
Suppose that there exist three constants $\e>0$ and $\alpha,N>1$
such that for every ball $B=B(x_0,r)$ with $\alpha B=B(x_0,\alpha r)
\subset E$, there exists a function 
$h=h_B\in L^p(B)$ with the properties:
$$
\align
\left\{\frac{1}{|B|}
\int_{B} |f-h|^2\, dx\right\}^{1/2}
& \le \e\,
\left\{\frac{1}{|\alpha B|}\int_{\alpha B}
|f|^2\, dx\right\}^{1/2},\tag 4.14\\
\left\{ \frac{1}{|B|}
\int_B |h|^p \, dx\right\}^{1/p}
&\le N \,
\left\{\frac{1}{|\alpha B|}
\int_{\alpha B} |f|^2 \, dx\right\}^{1/2}.
\tag 4.15
\endalign
$$
Then, if $2<q<p$ and $0<\e<\e_0=\e_0(n,p,q,\alpha,N)$, we have
$$
\left\{\frac{1}{|B|}\int_B |f|^q \, dx\right\}^{1/q}
\le C\, \left\{ \frac{1}{|\alpha B|}
\int_{\alpha B} |f|^2\, dx\right\}^{1/2},
\tag 4.16
$$
for any ball $B$ with
$\alpha B\subset E$,
where $C$ depends only on $n$, $p$, $q$, $\alpha$ and $N$.
\endproclaim

We remark that Theorem 4.13, whose proof is omitted
here,  is essentially proved in \cite{6}.

\demo{\bf Proof of Theorem C}
It suffices to show that $\CL$ satisfies condition (i)
in Theorem B for some $p>3$ if $n\ge 3$, and $p>4$ if $n=2$.
 This will be done by
combining Lemma 4.7 with Theorem 4.13.

Let $p$ and $\phi(r)$ be the same as in Lemma 4.7.
Fix $q$ so that $3<q<p$ if $n\ge 3$, and $4<q<p$ for $n=2$.
Choose $r_0>0$ so small that $\sup_{0<r<r_0} \phi(r)
<\e_0$, where $\e_0=\e_0(p,q,n,\alpha,C)$ is given in Theorem 4.13.

Now let $u$ be a weak solution of $\CL u=0$ in $D(x_0,\alpha_2 r)$
and $u=0$ on $B(x_0,\alpha_2r)\cap\partial\OO$, where $0<r<r_0/\alpha_2$.
If $B(x_0,\alpha_2 r)\subset \OO$, we simply apply Theorem 4.13
to $f=\nabla u$ on $E=B(x_0,\alpha_1 r)$ with
$h=\nabla v$ for each ball in $E$, given in Lemma 4.7.

In the case $x_0\in\partial\OO$, we also take $E=B(x_0,\alpha_1 r)$;
 but extend $f=\nabla u$ to be zero outside of $\OO$.
Given any $B^\prime=B(y_0,t)$ with $\alpha B^\prime \subset E$.
If $y_0\in \overline{\OO}$, we let $h_{B^\prime}
=\nabla v$ on $D(y_0,t)$
and zero otherwise.
If $y_0\notin \overline{\OO}$ and $B^\prime\cap\OO
\neq\emptyset$, we may find a ball $\widetilde{B}=B(z_0,2t)$
centered on $\partial\OO$,
 such that $B^\prime\subset \widetilde{B}$.
We  let $h_{B^\prime}=\nabla v$ on $D(z_0,2t)$ and zero
otherwise. Thus
the desired estimates (4.14)-(4.15)
for $B^\prime$ follow from (4.8)-(4.9).
This completes the proof.
\enddemo

\bigskip

\centerline{\bf 5. Proof of Theorems 3.1 and 3.2}

In this section we give the proof of Theorems 3.1 and 3.3,
 using a line
of argument similar to that in \cite{6}. 

\demo{\bf Proof of Theorem 3.1} 
Let $T$ be a bounded sublinear
operator
on $L^2(\br^n)$.
Suppose that $T$  satisfies assumption (3.2). 
We first note that
 with possibly different constants $\alpha_1,\alpha_2,N$, 
one may change balls $B$ in (3.2)
to cubes $Q$.

Fix $q\in (2,p)$.
Let $f$ be a bounded measurable function with
compact support. For $\lambda>0$, we consider the set
$$
E(\lambda)=\left\{ x\in\br^n: M(|Tf|^2)(x)>\lambda\right\},
\tag 5.1
$$
where $M$ is the Hardy-Littlewood maximal operator
defined by using cubes.
Since $Tf\in L^2$, $|E(\lambda)|\le C\, \| Tf\|_2^2/\lambda <\infty$. 
Let $A=1/(2\delta^{2/q})>5^n$, where
$\delta\in (0,1)$ is a small constant to be 
determined.
 Applying the Calder\'on-Zygmund
decomposition to $E(A\lambda)$, we obtain a collection
of disjoint dyadic cubes $\{ Q_k\}$ with the following
properties:
{(a)} $|E(A\lambda)\setminus \cup_k Q_k|=0$,\ \ 
{(b)} $|E(A\lambda)\cap Q_k|> \delta |Q_k|$,\ \ 
{(c)} $|E(A\lambda)\cap \overline{Q}_k |\le \delta|\overline{Q}_k|$,
where $\overline{Q}_k$ denotes the dyadic ``parent''of $Q_k$, i.e.,
$Q_k$ is one of the $2^n$
cubes obtained by bisecting the sides of $\overline{Q}_k$.
To see this, we first choose a large grid of dyadic cubes
of $\br^n$ so that $|Q\cap E(A\lambda)|<\delta |Q|$ for
each $Q$ in the grid. We then proceed as in the proof of
Lemma 1.1 in \cite{CP} for each $Q\cap E(A\lambda)\subset Q$.

We claim  that it is possible to choose constants $\delta, \gamma>0$
so that
$$
|E(A\lambda)|
\le \delta|E(\lambda)|
+|\left\{ 
x\in\br^n: \, M(|f|^2)(x)>\gamma\lambda\right\}|
\ \ \text{ for any }\lambda >0.
\tag 5.2
$$
This would imply that for any $\lambda_0>0$,
$$
\int_0^{A\lambda_0} \lambda^{\frac{q}{2}-1} |E(\lambda)|\, d\lambda
\le \delta A^{q/2}
\int_0^{\lambda_0}
\lambda^{\frac{q}{2}-1} |E(\lambda)|\, d\lambda
+C(\delta,\gamma)
\int_{\br^n} |f|^q\, dx.
\tag 5.3
$$
Using $\delta A^{q/2}=1/2^{q/2}<1$, $A>1$ and $\sup_{\lambda>0}
\lambda |E(\lambda)|<\infty$, we obtain
$$
\int_0^{\lambda_0}
\lambda^{\frac{q}{2}-1}
|\{ x\in \br^n:\,
|Tf(x)|^2>\lambda\}| \, d\lambda \le 
\int_0^{\lambda_0}
\lambda^{\frac{q}{2}-1} |E(\lambda)|\, d\lambda
\le C\, \int_{\br^n} |f|^q\, dx.
\tag 5.4
$$
Letting $\lambda_0\to\infty$ in (5.4), we conclude that
$\| Tf\|_q\le C\, \| f\|_q$.

It remains to prove (5.2). To this end, it suffices to show that
it is possible to choose
$\delta, \gamma>0$ such that if 
$\overline{Q}_k\cap \left\{x\in\br^n:  M(|f|^2)(x)
\le \gamma \lambda\right\}
\neq \emptyset$, then $\overline{Q}_k
\subset E(\lambda)$.
For this would imply that
$$
\aligned
|E(A\lambda)\cap \big\{ x\in \br^n:
&  M(|f|^2)(x) \le \gamma\lambda\big\}|
\le \sum_{k^\prime}
|E(A\lambda)\cap \overline{Q}_{k^\prime}|\\
&\le
\delta\sum_{k^\prime}
|\overline{Q}_{k^\prime}|
\le \delta |E(\lambda)|,
\endaligned
\tag 5.5
$$
where $\{ \overline{Q}_{k^\prime}\}$ is a disjoint subcover
of $E(A\lambda)\cap \big\{ x\in \br^n:\, M(|f|^2)(x)
\le \gamma
\lambda\big\}$ with the property that
$\overline{Q}_{k^\prime}\cap
\big\{ x\in\br^n: \,  M(|f|^2)(x) \le \gamma \lambda\big\}
\neq\emptyset$.

To finish the proof, we proceed by contradiction.
Suppose that there exists $x_0\in \overline{Q}_k
\setminus E(\lambda)$ and
$\{ x\in \overline{Q}_k :\, M(|f|^2)(x)
\le \gamma \lambda\}
\neq \emptyset$.
Then, if $Q$ contains $\overline{Q}_k$, we must have
$$
 \frac{1}{|Q|}\int_Q |f|^2 dx
\le \gamma \lambda\ \text{ and }\ 
\frac{1}{|Q|}\int_Q |Tf|^2 dx\le \lambda.
\tag 5.6
$$
It follows that for $x\in Q_k$,
$$
M(|Tf|^2)(x)
\le \max (M_{2\overline{Q}_k}(|Tf|^2)(x),
5^n\lambda),
\tag 5.7
$$
where $M_Q $ is a localized maximal function defined by
$$
M_Q (g)(x)
=\sup\Sb Q^\prime \owns x\\ Q^\prime\subset Q
\endSb
\frac{1}{|Q^\prime|}
\int_{Q^\prime} |g(y)| dy\ \ \ \text{ for }x\in Q.
\tag 5.8
$$
Since $A=1/(2\delta^{2/q})
\ge 5^n$, we have 
$$
\aligned
|Q_k\cap E(A\lambda)|
&\le |\big\{ x\in Q_k: M_{2\overline{Q}_k}(|Tf|^2)(x)
>A\lambda\big\} | \\
&
\le \big|\left\{ x\in Q_k: M_{2\overline{Q}_k}
( |T(f\chi_{
\alpha_2  \overline{Q}_k})|^2) (x)>\frac{A\lambda}{4}\right\}\big|\\
&\ \ \ \ \ \ \ 
+\big|\left\{
x\in Q_k : M_{2\overline{Q}_k}
(|T(f\chi_{\br^n\setminus \alpha_2 \overline{Q}_k})|^2)(x)>\frac{A\lambda}{4}
\right\} \big|\\
&
\le \frac{C_n}{A\lambda}
\int_{2\overline{Q}_k}
|T(f\chi_{\alpha_2 \overline{Q}_k})|^2 dx
+\frac{C_{n,p}}{(A\lambda)^{p/2}}
\int_{2\overline{Q}_k}
|T(f\chi_{\br^n\setminus \alpha_2 \overline{Q}_k})|^p dx.
\endaligned
\tag 5.9
$$
It then follows from the $L^2$ boundedness of $T$, assumption
(3.2) and (5.6)  that for any $\lambda>0$,
$$
\aligned
|Q_k\cap E(A\lambda)|
&\le |Q_k|
\left\{ \frac{C\gamma}{A} +\frac{C}{A^{p/2}}\right\}\\
&=\delta |Q_k|
\left\{ 2C\gamma \delta^{\frac{2}{q}-1}
+C2^{p/2} \delta^{\frac{p}{q}-1}\right\}.
\endaligned
\tag 5.10
$$
where $C$ depends only on $n,p,\alpha_1,\alpha_2,N$ as well as
 the operator norm of
$T$ on $L^2(\br^n)$.

Finally we choose $\delta\in (0,1)$ so small that
$C 2^{p/2} \delta^{\frac{p}{q}-1}\le 1/2$ and $A=1/(2\delta^{2/q})
\ge 5^n$.
This is possible since $q<p$.
With $\delta$ fixed, we choose $\gamma>0$ so small that 
$2C\gamma \delta^{\frac{2}{q}-1}\le 1/2$.
It follows from (5.10) that
$|Q_k \cap E(A\lambda)|\le \delta |Q_k|$. 
This contradicts with the fact that 
$|Q_k\cap E(A\lambda)|>\delta |Q_k|$.
Thus we must have $\overline{Q}_k
\subset E(\lambda)$ whenever the set $\{ x\in \overline{Q}_k:
M(|f|^2)(x)\le \gamma\lambda\}$ is not empty.
The proof is complete.
\enddemo

\remark{\bf Remark 5.11}
Let $T$ be a linear operator with kernel $K(x,y)$
satisfying
$$
|K(x,y)-K(x+h,y)|\le \frac{C|h|^\eta}{|x-y|^{n+\eta}},
\tag 5.12
$$
where $x,y,h\in \br^n$ and $|h|<|x-y|/4$.
Suppose supp$f\subset \br^n\setminus 8B$. Then
$$
|Tf(x)-Tf(y)|\le C\, \sup_{Q^\prime\supset Q}
\frac{1}{|Q^\prime|}
\int_{Q^\prime}
|f(z)| dz \ \ \ \text{ for any } x,y \in Q. 
\tag 5.13
$$
It follows that
$$
\| Tf\|_{L^\infty(Q)}
\le \frac{1}{|Q|}
\int_{Q} |Tf| dx
+C\, \sup_{Q^\prime\supset Q}
\frac{1}{|Q^\prime|} \int_{Q^\prime}|f| dx.
\tag 5.14
$$
Thus $T$ satisfies assumption (3.2) in Theorem 3.1 for any
$p>2$. Consequently, if $T$ is bounded on $L^2$,
then it is bounded on $L^p$
for any $2<p<\infty$. In this regard, Theorem 3.1 may be
considered as an extension of the Calder\'on-Zygmund Lemma.
\endremark

The following is a weighted version of Theorem 3.1.
Its proof may be carried out by a careful inspection of the
proof of Theorem 3.1. The key observation is that if
$d\mu =\oo^\delta dx$ where $\oo\in A_1(\br^n)$ and $0<\delta<1$, then
$\mu(E)\le C\left(\frac{|E|}{|Q|}\right)^{1-\delta}
\mu(Q)$ whenever $E\subset Q$.
We leave the details to the reader.

\proclaim{Theorem 5.15}
Under the same assumption as in Theorem 3.1, 
 $T$ is bounded on
$L^2(\br^n, \oo^\delta dx)$
where $\oo\in A_1(\br^n)$ and $0<\delta<1-\frac{2}{p}$.
\endproclaim

\demo{\bf Proof of Theorem 3.3}
The proof is similar to that of Theorem 3.1.
We first note that with possibly
different constants $\alpha_1,\alpha_2,N,r_0$,
inequality (3.4) holds for any ball $B(x_0,r)$
with the property that $0<r<r_0$ and $B(x_0,r)\cap\OO\neq\emptyset$.
Also one may replace balls $B$ in (3.4) by cubes
$Q$ of side length $r$.

Next we choose a cube $Q_0$ such that $\OO\subset Q_0$.
Fix $q\in (2,p)$. Let $\delta\in (0,1)$ be a small constant
to be determined. For $\lambda>0$, we consider the set
$$
E(\lambda)=\big\{ x\in Q_0:\
M_{2Q_0}(|Tf|^2\chi_\OO)(x)>\lambda\big\}.
\tag 5.16
$$
Then $|E(\lambda)|\le C\| f\|_2^2/\lambda \le \delta |Q_0|$ if
$$
\lambda \ge \lambda_1 =\frac{C}{\delta |Q_0|}\int_\OO |f|^2\, dx.
\tag 5.17
$$
Let $A=1/(2\delta^{2/q})$. For $\lambda\ge \lambda_1$,
we apply the Calder\'on-Zygmund
decomposition to $E(A\lambda)$. This produces a collection of
dyadic subcubes $\{ Q_k\}$ of $Q_0$ satisfying the same
properties (a), (b), (c) as in the proof of Theorem 3.1.
Note that $\delta |Q_k|<|E(A\lambda)|\le C\| f\|_2^2/(A\lambda)
\le \delta |Q_0|/A$. It follows that $|Q_k|< |Q_0|/A$.
 Thus we may choose  $\delta$ so small that
the side length of $2\overline{Q}_k$ is less than $r_0$.
With this observation, we may use the same argument as in the proof
of Theorem 3.1 to show that for $\lambda\ge \lambda_1$,
$$
|E(A\lambda)|\le \delta |E(\lambda)|
+|\big\{ x\in \br^n:\
M(|f|^2\chi_\OO)(x)>\gamma \lambda\big\}|.
\tag 5.18
$$
By integration, this implies that
$$
\int_\OO |Tf|^q \, dx
\le C\, \lambda_1^{q/2} |Q_0|
+C\, \int_\OO |f|^q \, dx
\le C\, \int_\OO |f|^q\, dx.
\tag 5.19
$$
The proof is finished.
\enddemo


\Refs
\widestnumber\key{19}

\ref\key 1
\by P.~Auscher
\paper
On necessary and sufficient conditions
for $L^p$ estimates of Riesz transform associated
to elliptic operators on $\br^n$ and related estimates
\jour Preprint
\yr 2004
\endref
 
\ref\key 2
\by P.~Aucher, T.~Coulhon, X.T.~Duong, and S.~Hofmann,
\paper Riesz transforms on manifolds and heat kernel regularity
\jour Preprint
\yr 2003
\endref


\ref\key 3
\by P.~Auscher and M.~Qafsaoui
\paper
Observation on $W^{1,p}$ estimates for divergence
elliptic equations with VMO coefficients
\jour Boll. Unione Mat. Ital. Sez. B Artic. Ric. Mat. (7)
\vol 5
\yr 2002
\pages 487-509
\endref

\ref\key 4
\by P.~Auscher and Ph.~Tchamitchian
\book Square Root Problem for Divergence Operators
and Related Topics
\bookinfo Ast\'erisque {\bf 249}
\publ Soc. Math. France
\yr 1998
\endref

\ref\key 5
\by P.~Auscher and Ph.~Tchamitchian
\paper Square roots of elliptic second order
divergence operators on strongly Lipschitz domains:
$L^p$ theory
\jour Math. Ann \vol 320 \pages 577-623
\yr 2001
\endref


\ref\key 6
\by L.A.~Caffarelli and I.~Peral
\paper On $W^{1,p}$ estimates for elliptic equations
in divergence form
\jour Comm. Pure App. Math. \vol 51
\yr 1998
\pages 1-21
\endref

\ref\key 7
\by T.~Coulhon and X.T.~Duong
\paper Riesz transforms for $1\le p\le 2$
\jour Trans. Amer. Math. Soc.
\vol 351 
\pages 1151-1169
\yr 1999
\endref

\ref \key 8
\by B.~Dahlberg and C.~Kenig
\paper Hardy spaces and the Neumann problem in $L^p$
for Laplace's equation in Lipschitz domains
\jour Ann. of Math. \vol 125\yr 1987
\pages 437-466
\endref

\ref\key 9
\by E.~Davies
\book Heat Kernels and Spectral Theory
\publ Cambridge Univ. Press
\yr 1989
\endref

\ref\key 10
\by J.~Duoandikoetxea
\book Fourier Analysis
\bookinfo Graduate Studies in Math.
\vol 29 \publ Amer. Math. Soc.
\yr 2000
\endref

\ref\key 11
\by G.~Di Fazio
\paper $L^p$ estimates for divergence form elliptic
equations with discontinuous coefficients
\jour Boll. Un. Mat. Ital. A (7)
\vol 10 \yr 1996
\pages 409-420
\endref

\ref\key 12
\by F. Giaquinta
\book Multiple Integrals in the Calculus of Variations and
Non-Linear Elliptic Systems
\bookinfo Annals of Math. Studies
\vol 105
\publ Princeton Univ. Press
\yr 1983
\endref

\ref\key 13
\by D.~Jerison and C.~Kenig
\paper
The inhomogeneous Dirichlet problem in Lipschitz domains
\jour J. Funct. Anal. \vol 130
\pages 161-219
\yr 1995
\endref

\ref\key 14
\by C.~Kenig
\book Harmonic Analysis Techniques for
Second Order Elliptic Boundary Value Problems
\bookinfo Regional Conf. Series in Math., no. 83
\publ Amer. Math. Soc.
\yr 1994
\endref

\ref \key 15
\by N.~Meyers
\paper An $L^p$ estimate for the gradient of solutions of second
order elliptic divergence equations
\jour Ann. Sc. Norm. Sup. Pisa
\vol 17 \yr 1963 \pages 189-206
\endref

\ref \key 16
\by J.L.~Rubio de Francia
\paper Factorization theory and the $A_p$ weights
\jour Amer. J. Math.
\vol 106 \yr 1984 \pages 533-547
\endref

\ref \key 17
\by Z.~Shen
\paper The $L^p$ Dirichlet problem for the elliptic system
on Lipschitz domains
\jour Preprint
\yr 2004
\endref

\ref\key 18
\by E.~Stein
\book Singular Integrals and Differentiability
Properties of Functions
\publ Princeton Univ. Press
\yr 1970
\endref

\ref\key 19
\by L.~Wang
\paper A geometric approach to the Calder\'on-Zygmund estimates
\jour Acta Math. Sinica (Engl. Ser.) \vol 19
\yr 2003 \pages 381-396
\endref

\endRefs
\enddocument

\end